\documentclass[11pt]{article}
\usepackage[english]{babel}
\usepackage{amsthm,amsmath,amssymb,amsfonts,latexsym}
\newtheorem{theorem}{\qquad Theorem}
\newtheorem{corollary}{\qquad Corollary}
\newtheorem{lemma}{\qquad Lemma}

\title{Solutions of polynomial equation over $\mathbb{F}_p$ and new bounds of additive energy\footnote{The research was carried out at the IITP RAS at the expense of the Russian
Foundation for Sciences (project N~14-50-00150).}}

\author{Ilya Vyugin and Sergey Makarychev}

\date{}

\begin{document}
\maketitle

\begin{abstract}
We present a new proof of  Corvaja and Zannier's \cite{C-Z} the upper bound of the number of solutions $(x,y)$ of the algebraic equation $P(x,y)=0$  over a field $\mathbb{F}_p$ ($p$ is a prime), in the case, where $x\in g_1G$, $y\in g_2G$, ($g_1G$, $g_2G$ -- are cosets by some subgroup $G$ of a multiplicative group $\mathbb{F}_p^*$). The estimate of  Corvaja and Zannier was improved in average,
and some applications of it has been obtained. In particular we present the new bounds of additive and polynomial energy. 
\end{abstract}

\section{Introduction}

We study an algebraic equation
\begin{eqnarray}\label{equa-first}
P(x,y)=0
\end{eqnarray}
over a field $\mathbb{F}_p$ (or its algebraic closure $\overline{\mathbb{F}}_p$), where $p$ is a prime number. Suppose that $P\in\mathbb{F}_p[x,y]$ is an absolutely irreducible polynomial of two variables $x$ and $y$. 
Let $G$ be a subgroup of $\mathbb{F}_p^*$ (multiplicative group of $\mathbb{F}_p$). We study the upper bound of the number  solutions of equation (\ref{equa-first}), such that $x\in g_1G$, $y\in g_2G$.


The first result of such a type belongs to Garcia and Voloch \cite{GV-1}. Their result has been improved by Heath-Brown and Konyagin \cite{HK-1}. They proved using Stepanov method (see \cite{HK-1},\cite{Step}) that for any subgroup $G\subset \mathbb{F}_p^*$, such that $|G|<(p-1)/((p-1)^{1/4}+1)$ and an arbitrary nonzero $\mu$ the number of solutions of the linear equation 
\begin{eqnarray}\label{equa-linear}
y=x+\mu
\end{eqnarray}
such that $(x,y)\in G\times G$, does not exceed $4|G|^{2/3}$. In the other words they studied such a problem for linear equations (\ref{equa-first}). The case of such systems was studied in \cite{VS},\cite{SSV}.

Corvaja and Zannier \cite{C-Z} have obtained the following theorem.
\begin{theorem}[Corvaja and Zannier, \cite{C-Z}]\label{CZ-th}
Let $X$ be a smooth projective absolutely irreducible curve over a field $\kappa$ of characteristic $p$. Let $u,v\in \kappa(X)$ be rational functions, multiplicatively independent modulo $\kappa^{*}$, and with non-zero differentials; let $S$ be the set of their zeros and poles; and let $\chi=|S|+2g-2$ be the Euler characteristic of $X\setminus S$. Then
\begin{eqnarray}\label{poly-form}
\sum_{\nu\in X(\overline{\kappa})\setminus S}\min\{ \nu(1-u),\nu(1-v)\}\leqslant\left(3\sqrt[3]{2}(\deg u \deg v)^{1/3},12\frac{\deg u \deg v}{p}\right),
\end{eqnarray}
where $\nu(f)$ denotes the multiplicity of vanishing of $f$ at the point $\nu$.
\end{theorem}

Corollary 2 of the paper \cite{C-Z} gives us the estimate
$$
\#|\{ (x,y)\mid P(x,y)=0, \, (x,y)\in g_1G\times g_2G\}|\leqslant \max\left( 3\sqrt[3]{2}(mn\chi)^{1/3}|G|^{2/3},12\frac{mn|G|^{2/3}}{p}\right),
$$
where $(m,n)$ is a bidegree of the polynomial $P(x,y)$, $\chi$ is the Euler characteristic of the curve (\ref{equa-first}), $g_1G$ and $g_2G$ are costes by a subgroup $G$.
The proof of Theorem \ref{CZ-th} using Wronskians method over fields of positive characteristics. At first we give a new proof of such a type estimate using Stepanov method (see Theorem \ref{th1}) and at the second we modify the proof of Theorem \ref{CZ-th} and improve the estimate in average by Konyagin's modification of Stepanov method (see Theorem \ref{th-sr}).

\begin{theorem}\label{th1}
Let $P(x,y)$ be a polynomial of the form (\ref{equa-first})), of bidegree $(m,n)$, $P(0,0)\not=0$, $\deg P(x,0)\geqslant 1$, $G$  is a subgroup of $\mathbb{F}_p^*$, $100(mn)^{3/2}<|G|<\frac{1}{3}p^{3/4}$, $g_1,g_2\in \mathbb{F}_p^{*}$, and
\begin{eqnarray}\label{M2}
M_1=\{ (x,y)\mid P(x,y)=0,\,\, x\in g_1G, \,\,  y \in g_2G  \}.
\end{eqnarray}
Then the following estimate holds $\# M_1\leqslant 16mn^2(m+n)|G|^{2/3}$.
\end{theorem}

The new results of the paper concerns with the estimates in average and estimates of the additive and polynomial energies. Let us consider a homogeneous polynomial $P(x,y)$ of degree $n$. We estimate a sum $L$ of numbers of solutions of the set of equations:
\begin{eqnarray}\label{set-equ}
P(x,y)=l_i,\qquad i=1,\ldots, h.
\end{eqnarray}
 Theorem \ref{th1} gives us the trivial bound $L\leqslant 16mn^2h(m+n)|G|^{2/3}$. 


\begin{theorem}\label{th-sr}
Let us consider a homogeneous polynomial $P(x,y)$ of degree $n$, such that $\deg P(x,0)\geqslant 1$, $P(0,0)\not=0$ and a set of equations (\ref{set-equ}) such that $l_1,\dots,l_h$ belong to different cosets $g_iG$, and $h<\min(\frac{1}{81}|G|^{4/3},\frac{1}{3}pt^{-4/3})$. Then the sum $N_h$ of numbers of solutions of the set of equations (\ref{set-equ}) does not exceed $32h^{3/4}n^5|G|^{2/3}$.
\end{theorem}

\section{Corollaries and applications}





Let $A,B$ be subsets of a field $\mathbb{F}_p$. The {\it additive energy} is defined by
$$
E(A,B)=\#\{ (x_1,y_1,x_2,y_2)\mid x_1+y_1=x_2+y_2,\,\,\, x_1,x_2\in A,\,\, y_1,y_2\in B\},
$$
and we denote $E(A,A)$ by $E(A)$. The additive energy plays an important role in many problems of additive combinatorics as well as in number theory (see e.g. \cite{TV},\cite{ss_E_k}).

We consider some generalization of the additive energy which we call a {\it polynomial energy}. Polynomial energy is the following 
$$
E_P^q(A)=\#\{ (x_1,y_2,x_2,y_2)\mid P(x_1,y_1)=P(x_2,y_2),\,\,\, x_1,y_1,x_2,y_2\in A\},
$$
where $P(x,y)\in \mathbb{F}_p[x,y]$ is a polynomial. We will consider polynomials $P(x,y)$ of bidegree $(m,n)$ such that $\deg P(x,0)\geqslant 1$.

\begin{corollary}
Let $P(x,y)$ be a polynomial of bidegree $(m,n)$ such that $\deg P(x,0)\geqslant 1$ and $G$ be a subgroup of $\mathbb{F}_p^*$. Then the number of solutions $(x,y,z,w)$ of the equation
$$ 
P(x,y)=P(z,w)
$$
such that $x,y,z,w\in G$,
does not exceed $17mn^2(m+n)|G|^{8/3}$.
\end{corollary}

{\it Proof.}
Let us fix two variables, for example, $z$ and $w$. Theorem \ref{th1} gives us that if $P(0,0)-P(z,w)\not=0$, then the number of solutions $(x,y)$ of the equation $P(x,y)=P(z,w)$ does not exceed $16mn^2(m+n)|G|^{\frac{2}{3}}$. The condition  $P(0,0)-P(z,w)\not=0$ can be not satisfied only for $n|G|$ pairs $(z,w)\in G\times G$.
Note that for each fixed $z$ and $t$ the number of solutions does not exceed $16mn^2(m+n)|G|^{2/3}$ if $P(0,0)-P(z,w)\not=0$. So let us obtain that  the number of solutions of the polynomial equation 
$$ 
P(x,y)=P(z,t)
$$
does not exceed $16mn^2(m+n)|G||G||G|^{2/3}+n^2|G|^2\leqslant 17mn^2(m+n)|G|^{8/3}$.
$\Box$

\begin{theorem}\label{th-energy}
Let us suppose that  $100(mn)^{3/2}<|G|<\left(\frac{p}{3}\right)^{\frac{12}{17}}$. Then the following holds:
if $q\leqslant 3$ then
$$
E_P^q(G)\leqslant C(n,q)|G|^{\frac{7q+16}{12}};
$$
if $q=4$ then 
$$
E_P^4\leqslant C(n,q)|G|^{1+\frac{2q}{3}}\ln |G|;
$$
if  $q\geqslant 5$ then
$$
E_P^q(G)\leqslant C(n,q)|G|^{1+\frac{2q}{3}},
$$
where $C(n,q)$ depends only on $n$ and $q$.
\end{theorem}


Let us consider the sets $f(G)=\{ f(x)\mid x\in G\}$ and $g(G)=\{ g(x)\mid x\in G\}$, where $G$ is a subgroup of $\mathbb{F}_p^*$, and $f,g\in\mathbb{F}_p[x]$.

\begin{corollary}
Let $G$ be a subgroup of $\mathbb{F}_p^*$, and $f,g\in\mathbb{F}_p[x]$, $\deg f=m$, $\deg g=n$ and  $100(mn)^{3/2}<|G|<\frac{1}{3}p^{3/4}$. Then
$$
E(f(G),g(G))\leqslant 16mn^2(m+n)|G|^{8/3}.
$$
\end{corollary}

{\it Proof.} It is easy to see that 
\begin{eqnarray}
E(f(G),g(G))\leqslant |\{ (x_1,y_1,x_2,y_2)\mid f(x_1)+g(y_1)=f(x_2)+g(y_2),\,\, x_1,x_2,y_1,y_2\in G \}|=\\
=|\{ (x_1,y_1,x_2,y_2)\mid f(x_1)-f(x_2)=\mu=g(y_2)-g(y_1),\,\, x_1,x_2,y_1,y_2\in G, \mu\in \mathbb{F}_p \}|.\label{set2}
\end{eqnarray}
We obtain the following
$$
E(f(G),g(G))\leqslant 16mn^2(m+n)^2|G|^{2/3}|G|^2=16mn^2(m+n)|G|^{8/3}.
$$
$\Box$





\section{Proof of Theorem \ref{th1}}

\subsection{Stepanov method with polynomials of two variables}

Let us consider a polynomial $\Phi\in \mathbb{F}_p[X,Y,Z]$ such that
$$
\deg_X\Phi<A,\quad \deg_Y\Phi<B,\quad \deg_Z\Phi<C,
$$
or in the other words
\begin{eqnarray}\label{poly-step2}
\Phi(X,Y,Z)=\sum_{a,b,c}\lambda_{a,b,c}X^aY^bZ^c,\qquad a\in[A],\quad b\in[B],\quad c\in[C],
\end{eqnarray}
where $[N]=\{ 0,1,\ldots, N-1\}$.
Consider the following polynomial
\begin{eqnarray}\label{Psi2}
\Psi(x,y)=\Phi(x,x^t, y^t),
\end{eqnarray}
which satisfies to the following conditions:

1) all roots $(x,y)$, such that $x\in g_1G,y\in g_2G$, of the equation (\ref{equa-first}) are zeros of system
\begin{eqnarray}\label{Syst}
\left\{
\begin{aligned}
&\Psi(x,y)=0 \\
&P(x,y)=0\\
\end{aligned}
\right.
\end{eqnarray}
of the order at least $D$.

2) the greatest common divisor of polynomials $\Psi(x,y)$ and $P(x,y)$ is a constant. 

If these conditions are satisfied then the generalized B\'ezout's theorem gives us the upper bound of the number $N$ of roots $(x,y)$ such that $x\in g_1G$, $y\in g_2G$:
\begin{eqnarray}\label{ord-N}
N\leqslant \frac{\deg \Psi(x,y) \cdot\deg P(x,y)}{D}\leqslant  \frac{(A-1+(B-1)t+(C-1)t)(m+n)}{D}.
\end{eqnarray}

A pair $(x,y)$ is the solution of the system (\ref{Syst}) of the order at least $D$, if $P(x,y)=0$ and $\Psi(x,y)=0$ and derivatives
$$
\frac{d^k}{dx^k}\Psi(x,y)=0,\qquad k=1,\ldots,D-1
$$
vanishes on the curve $P(x,y)=0$.

\subsection{Lemmas}

\begin{lemma}\label{lem-div}
Let $Q(x,y)\in \mathbb{F}_p[x,y]$ be a polynomial such that
$$
\deg_x Q(x,y)\leqslant \mu,\quad \deg_y Q(x,y)\leqslant \nu
$$
and $P(x,y)\in \mathbb{F}_p[x,y]$ be a polynomial such that
$$
\deg_x P(x,y)\leqslant m,\quad \deg_y P(x,y)\leqslant n.
$$
Then the condition
$$
P(x,y)\mid Q(x,y)
$$
on coefficients of the polynomial $Q(x,y)$ can be given by $n((\nu-n+2)m+\mu)\leqslant (\mu+\nu+1)mn$ homogeneous linear algebraic equations.
\end{lemma}

{\it Proof.} Consider the polynomial
$$
P(x,y)=f_n(x)y^n+\ldots+f_1(x)y+f_0(x),\qquad \deg f_i(x)\leqslant m
$$
and the polynomial 
$$
Q_0(x,y)=Q(x,y)f_n(x)=g_{0,\nu}(x)y^{\nu}+\ldots+g_{0,1}(x)y+g_{0,0}(x).
$$ 
Let us construct the polynomials $Q_i(x,y)=g_{i,\nu-i}(x)y^{\nu-i}+\ldots+g_{i,1}(x)y+g_{i,0}(x)$, $i=1,\ldots,\nu-n+1$ such that
$$
Q_i(x,y)=Q_{i-1}(x,y)-\frac{g_{i-1,\nu-i+1}(x)}{f_n(x)}P(x,y).
$$
It is easy to see that $\deg_y Q_{i}(x,y)< \deg_y Q_{i-1}(x,y)$ and $\frac{g_{i-1,\nu-i+1}(x)}{f_n(x)}$ --- is a polynomial, because $f_n(x)\mid g_{i-1,\nu-i+1}(x)$ and $\deg g_{i,j}(x)\leqslant \mu+(i+1)m$.

Consequently, $P(x,y)\mid Q(x,y)$ if and only if $Q_{\nu-n+1}(x,y)\equiv 0$. The polynomial $Q_{\nu-n+1}(x,y)$ has $n((\mu+(\nu-n+2)m)$ coefficients which are homogeneous linear forms of coefficients of polynomial $Q(x,y)$.
We have  $n((\nu-n+2)m+\mu)$  homogeneous linear algebraic equations.
$\Box$

\begin{lemma}\label{lem-div-nonzero}
Let us consider a polynomial $Q\in\mathbb{F}_p[x,y]$ and an irreducible polynomial 
$$
P(x,y)=f_n(x)y^n+\ldots+f_1(x)y+f_0(x)
$$ 
of bidegree $(m,n)$. If $P(x,y)\mid Q(x,y^t)$, and $t\mid (p-1)$, then $P(x,0)^{\lfloor t/n\rfloor}\mid Q(x,0)$.\footnote{$\lfloor x\rfloor$ -- is the largest integer less than $x$.}
\end{lemma}

{\it Proof.} We have the following $P(x,y)\mid Q(x,y^t)$. Let us substitute $y=g\tilde{y}$, where $g\in G$ -- is the group of $t$-roots of $1$, to the polynomial $P(x,y)\longmapsto P_g(x,\tilde{y})=P(x,g\tilde{y})$. It is easy to see that
$$
P_g(x,y)\mid Q(x,y^t)
$$
for any $g\in G$ polynomials $P_g(x,y)$ are irreducible. The leading coefficient of the polynomial $P_g(x,y)$ is $f_0(x)g^n$. There are at least $\lfloor t/n \rfloor$ elements $g_1,\ldots,g_{\lfloor t/n\rfloor}\in G$ such that $g_1^n,\ldots,g_{\lfloor t/n\rfloor}^n$ are pairwise distinct. Note that the free terms of the polynomials $P_{g_1}(x,y),\ldots, P_{g_{\lfloor t/n\rfloor}}(x,y)$ are the same.
Consequently, we have that the polynomials $P_{g_1}(x,y),\ldots, P_{g_{\lfloor t/n\rfloor}}(x,y)$ relatively prime, and
$$
(P_{g_1}(x,y)\ldots P_{g_{\lfloor t/n\rfloor}}(x,y))\mid Q(x,y^t).
$$
Here we have that
$$
(P_{g_1}(x,0)\ldots P_{g_{\lfloor t/n\rfloor}}(x,0))\mid Q(x,0),
$$
and considering that $P(x,0)=P_g(x,0)$ for any $g\in\mathbb{F}_p^*$ we obtain the statement of the Lemma 
$$
P(x,0)^{\lfloor t/n\rfloor}\mid Q(x,0).
$$
$\Box$

\begin{lemma}\label{lem-nonzero}
Let
$$
\Psi(x,y)=\sum_{a,b,c}\lambda_{a,b,c}x^ax^{bt}y^{ct},\qquad a\in [A],\quad b\in[B],\quad c\in[C],
$$
be a polynomial, $nAB\leqslant t$, coefficients $\lambda_{a,b,c}$ do not vanish simultaneously, $P(x,y)$ be an irreducible polynomial and $\deg_y P(x,y)=n$, $P(0,0)\not=0$. Then  $P(x,y)$ does not divide $\Psi(x,y)$. 
\end{lemma}

{\it Proof.} Let us denote $c_{min}=\min_{a,b,c:\lambda_{a,b,c}\not=0} c$. Consider the polynomial $\Psi$ of the form 
$$
\Psi(x,y)=y^{c_{min}t}\tilde{\Psi}(x,y).
$$
It is easy to see that
$$
\tilde{\Psi}(x,y)=\sum_{a,b,c:c>c_{min}}\lambda_{a,b,c}x^ax^{bt}y^{(c-c_{min})t}+\sum_{a,b}\lambda_{a,b,c_{min}}x^ax^{bt},\quad a\in[A],\, b\in[B],\, c\in[C].
$$
So, if $P(x,y)\mid \Psi(x,y)$ then $P(x,y)\mid \tilde{\Psi}(x,y)$ and
$$
P(x,0)^{\lfloor t/n\rfloor}\mid\tilde{\Psi}(x,0)
$$ 
by Lemma \ref{lem-div-nonzero} 
and $\tilde{\Psi}(x,0)\not=0$.
It can not be true if $P(x,0)$ has at least one nonzero root and the number of members of polynomial $\Psi(x,0)$ does not exceed $t/n$ ($t\geqslant nAB$). 
$\Box$

\subsection{Derivatives and differential operators}

Let us express derivatives $\frac{d^k}{dx^k}y$ on the algebraic curve $P(x,y)=0$.
Consider the polynomials $q_k(x,y)$ and $r_{k}(x,y)$, $k\in\mathbb{N}$, which are defined by induction
$$
q_1(x,y)=-\frac{\partial}{\partial x}P(x,y),\qquad r_1(x,y)=\frac{\partial}{\partial y}P(x,y),
$$
and
$$
q_{k+1}(x,y)=\frac{\partial q_k}{\partial x}\left(\frac{\partial P}{\partial y}\right)^2-\frac{\partial q_k}{\partial y}\frac{\partial P}{\partial x}\frac{\partial P}{\partial y}-(2k-1)q_k(x,y)\frac{\partial^2 P}{\partial x\partial y}\frac{\partial P}{\partial y}+(2k-1)q_k(x,y)\frac{\partial^2 P}{\partial y^2}\frac{\partial P}{\partial x},
$$
$$
r_{k+1}(x,y)=r_k(x,y)\left(\frac{\partial P}{\partial y}\right)^2=\left(\frac{\partial P}{\partial y}\right)^{2k+1},\qquad k=\mathbb{N}.
$$
Derivatives of algebraic function $y(x)$ which are defined by the equation $P(x,y)=0$ have the following expressions $\frac{d^k}{dx^k}y=\frac{q_k(x,y)}{r_k(x,y)}$, $k\in\mathbb{N}$.
Actually, we have the following expressions
$$
\frac{d}{dx}y=\frac{q_1(x,y)}{r_1(x,y)}=-\frac{\frac{\partial}{\partial x}P(x,y)}{\frac{\partial}{\partial y}P(x,y)},
$$
$$
\frac{d^{k+1}}{dx^{k+1}}y=\frac{q_{k+1}(x,y)}{r_{k+1}(x,y)}=\frac{
\frac{\partial q_k}{\partial x}\left(\frac{\partial P}{\partial y}\right)^2-\frac{\partial q_k}{\partial y}\frac{\partial P}{\partial x}\frac{\partial P}{\partial y}-(2k-1)q_k(x,y)\frac{\partial^2 P}{\partial x\partial y}\frac{\partial P}{\partial y}+(2k-1)q_k(x,y)\frac{\partial^2 P}{\partial y^2}\frac{\partial P}{\partial x}}{r_k(x,y)\left(\frac{\partial P}{\partial y}\right)^2}.
$$

Let us obtain the following lemma.

\begin{lemma}\label{degs-poly}
Degrees of polynomials $q_k(x,y)$ and $r_k(x,y)$ satisfy to the following bounds
$$
\deg_{x}q_k(x,y)\leqslant (2k-1)m-k,\quad \deg_{y}q_k(x,y)\leqslant (2k-1)n-k+1,
$$
$$ 
\deg_{x}r_k(x,y)\leqslant (2k-1)m, \quad \deg_{y}r_k(x,y)\leqslant (2k-1)(n-1),\quad k\in\mathbb{N}.
$$ 
\end{lemma}

{\it Proof. } It is easy to see that $\deg_{x}q_1(x,y)\leqslant m-1$, $\deg_{y}q_1(x,y)\leqslant n$ and
$$
\deg_{x}q_k(x,y)\leqslant\deg_{x}q_{k-1}(x,y)+2m-1\leqslant  (2k-1)m-k,
$$
$$
\deg_{y}q_k(x,y)\leqslant\deg_{y}q_{k-1}(x,y)+2n-1\leqslant  (2k-1)n-k+1.
$$
For the polynomial $r_{k}(x,y)$ the statement is obvious. $\Box$


Let us define the differential operators
\begin{eqnarray}\label{diff-oper}
D_k=\left(\frac{\partial P}{\partial y}\right)^{2k-1}x^ky^k\frac{d^k}{dx^k},\qquad k=\mathbb{N}.
\end{eqnarray}
It is easy to see that the following relations holds
\begin{eqnarray}\label{diff-oper-rel}
D_k x^ax^{bt}y^{ct}=R_{k,a,b,c}(x,y)x^{a}x^{bt}y^{ct},\\
D_k\Psi(x,y)|_{x,y\in G}=R_k(x,y)|_{x,y\in G},\nonumber
\end{eqnarray}
with some polynomials $R_{k,a,b,c}(x,y)$ and $R_k(x,y)$.
Let us obtain the following Lemma \ref{lem-R}.

\begin{lemma}\label{lem-R}
Degrees of polynomials $R_{k,a,b,c}(x,y)$ and $R_k(x,y)$ satisfy to the following bounds 
$$
\deg_x R_{k,a,b,c}(x,y)\leqslant 2(2k-1)m\leqslant 4km\qquad \deg_y R_{k,a,b,c}(x,y)\leqslant 2(2k-1)(2n-1)+1\leqslant 4kn
$$
$$
\deg_x R_{k}(x,y)\leqslant A+4km\qquad \deg_y R_{k}(x,y)\leqslant 4kn.
$$
\end{lemma}

{\it Proof.} This follows easily from Lemma \ref{degs-poly} and formulas (\ref{diff-oper}),(\ref{diff-oper-rel}). $\Box$

Let us consider the system
\begin{eqnarray}\label{resul}
\left\{
\begin{aligned}
&P(x,y)=0 \\
&\frac{\partial P}{\partial y}(x,y)=0\\
\end{aligned}
\right..
\end{eqnarray}
Polynomials $P(x,y)$ and $\frac{\partial P}{\partial y}(x,y)$ are relatively prime, because $P(x,y)$ is irreducible. It means that B\'ezout's theorem gives us the bound $L\leqslant (m+n)(m+n-1)$, where $L$ is the number of roots of the system (\ref{resul}) (see \cite{CKW}).

We have the following lemma.
\begin{lemma}\label{addzeros}
If $\Psi(x,y)=0$ and $D_j\Psi(x,y)=0$, $j=1,\ldots,k-1$ then at least one of the following alternatives holds: either 


- $(x,y)$ is a root of the order at least $k$ the system (\ref{resul});

- $x=0$ or $y=0$ or $\frac{\partial P}{\partial y}(x,y)=0$.
\end{lemma}

{\it Proof.} This is a direct consequence of the formula (\ref{diff-oper}). $\Box$

\subsection{End of the proof of Theorem \ref{th1}}

Let us suppose that $P(x,y)$ is the absolutely irreducible polynomial.
Define the following parameters
$$
A=\left\lfloor\frac{t^{2/3}}{n}\right\rfloor,\quad B=C=\lfloor t^{1/3}\rfloor
$$
$$
D=\left\lfloor \frac{B^2}{4mn^2}\right\rfloor .
$$
Consider the polynomial (\ref{Psi2}) and the system (\ref{Syst}). The condition
\begin{eqnarray}\label{Deriv-2}
D_k\Psi(x,y)=0\quad \text{if $P(x,y)=0$ and $(x,y)\in g_1G\times g_2G$, $k=0,\ldots,D-1$}
\end{eqnarray}
can be calculated by means of Lemmas \ref{lem-R} and \ref{lem-div}. The condition (\ref{Deriv-2}) is equivalent to the set 
of
$$
mn\sum_{k=0}^{D-1}(4km+4kn+A+1)=(A+1)Dmn+2mn(m+n)D(D-1)\leqslant ADmn+2mn(m+n)D^2
$$
homogeneous linear algebraic equations of variables $\lambda_{a,b,c}$ . This system has a nonzero solution if the inequality holds
\begin{eqnarray}\label{number-of-equ}
2D^2mn(m+n)+DmnA< ABC.
\end{eqnarray}
The inequality (\ref{number-of-equ}) is the following
$$
2D^2mn(m+n)+DmnA<t^{4/3}\frac{m+n+2mn}{8mn^3} <\frac{1}{2}\frac{t^{4/3}}{n}=\left\lfloor\frac{t^{2/3}}{n}\right\rfloor \lfloor t^{1/3}\rfloor^2=ABC,
$$
and it is satisfied if $t>8t^{3/2}$.
The conditions of Lemma \ref{lem-nonzero} hold
$$
t\geqslant nAB=n\left\lfloor\frac{t^{2/3}}{n}\right\rfloor\lfloor t^{1/3}\rfloor,
$$
and conditions
$$
\deg \Psi(x,y)<A+Bt+Ct<p,\qquad \deg P(x,y)<m+n<p
$$
hold too. It is easy to see that by (\ref{ord-N}) and Lemma \ref{addzeros} we obtain the following bound
$$
N\leqslant\frac{(m+n)(A+Bt+Ct)}{D}\leqslant 16mn^2(m+n)t^{2/3},
$$
because $t>100(mn)^{3/2}$ and, consequently, $\left\lfloor\frac{B^2}{4mn^2}\right\rfloor>\frac{B^2}{4mn^2}-1>\frac{3}{4}\frac{t^{2/3}}{4mn^2}$. 

Consider the case of reducible polynomial $P(x,y)$ over the field $\overline{\mathbb{F}}_p$. Consider the polynomial $P(x,y)$ as a product of irreducible polynomials $P_i(x,y)$:
$$
P(x,y)=\prod_{i=1}^sP_i(x,y).
$$
Let us denote degrees as following $\deg_x P_i(x,y)=m_i$, $\deg_y P_i(x,y)=n_i$, and $m=\sum_{i=1}^sm_i$, $n=\sum_{i=1}^sn_i$.
The set $M_1\subseteq \cap_{i=1}^s M_{1,i}$, where 
$$
M_{1,i}=\{ (x,y)\mid P_i(x,y)=0,\,\, x\in g_1G, \,\,  y \in g_2G  \}.
$$
Consequently, we have the estimate
$$
\# M_1\leqslant \sum_{i=1}^s 16m_in_i^2(m_i+n_i)|G|^{2/3}\leqslant 16mn^2(m+n)|G|^{2/3}.
$$
Theorem \ref{th1} is proved. 
$\Box$

\section{Proof of Theorem \ref{th-sr}}

Let us consider the equations 
\begin{eqnarray}\label{poly-equ-l}
P(x,y)=l 
\end{eqnarray}
and the equation 
\begin{eqnarray}\label{poly-equ-gamma}
P(x,y)=\gamma,
\end{eqnarray}
with $l$ and $\gamma$ such that $l\in \gamma\Gamma$, where $\Gamma$ is a subgroup of $n$-powers of $\mathbb{F}_p^*$ and $l/\gamma\notin G$. 
Then $l=\gamma \mu^n$ for some $\mu\in\mathbb{F}_p^*$ and the equation $P(x,y)=l=\gamma \mu^n$ is equivalent to the equation $P\left(\frac{x}{\mu},\frac{y}{\mu}\right)=\gamma$. Actually, if $x,y\in G$ then $\frac{x}{\mu},\frac{y}{\mu}\in \mu^{-1}G$. Consequently, the number of solutions of equation (\ref{poly-equ-l}) with restriction $x,y\in G$ is equivalent to the equation (\ref{poly-equ-gamma}) with restriction $x,y\in \mu^{-1}G$. Also, note that $\mu^{-1}G\cap G=\emptyset$.

Let us consider some set of equations (\ref{set-equ}). It is easy to see that $l_i=\gamma \mu_i^n$, $i=1,\dots,h$ and $\mu_i^{-1}G\cap\mu_j^{-1}G=\emptyset$, $i\not= j$. The sum of numbers of solutions (\ref{set-equ}) is equal to the number of solutions of equation (\ref{poly-equ-gamma}) with restriction $(x,y)\in \bigcup_{i=1}^h(\mu_i^{-1}G\times \mu_i^{-1}G)$.

Let us update Stepanov method. To construct a polynomial (\ref{poly-step2}) such that all roots of equations (\ref{set-equ}) be roots of (\ref{poly-step2}) of orders at least $D$.

Take the following parameters
$$
A=\lfloor h^{-1/2}t^{2/3}\rfloor,\quad B=C=\lfloor h^{1/4}t^{1/3}\rfloor
$$
$$
D=\left\lfloor h^{-1/2}\frac{t^{2/3}}{4n^3}\right\rfloor.
$$

The condition
\begin{eqnarray*}\label{Deriv-3}
D_k\Psi(x,y)=0\quad \text{if $P(x,y)=0$ and $(x,y)\in \bigcup_{i=1}^h(\mu_i^{-1}G\times \mu_i^{-1}G)$}
\end{eqnarray*}
can be calculated by means of Lemmas \ref{lem-R} and \ref{lem-div}. The condition (\ref{Deriv-2}) is equivalent to the set 
of homogeneous equations 
$$
n^2\sum_{k=0}^{D-1}(8kn+A+1)=(A+1)Dn^2+2n^3D(D-1)\leqslant ADn^2+2n^3D^2
$$
This system has a nonzero solution if
\begin{eqnarray}\label{number-of-equ-2}
h(ADn^2+4n^3D^2)< ABC.
\end{eqnarray}

Obviously, if $t>8h^{3/2}$ then the following holds 
\begin{eqnarray}\label{number-of-equ-3}
h\left(h^{-1}\frac{t^{2/3}}{4n^3}n^2+ h^{-1}4n^3\frac{t^{4/3}}{16n^6}\right)< \lfloor h^{-1/2}t^{2/3}\rfloor \lfloor h^{1/4}t^{1/3}\rfloor^2.
\end{eqnarray}
Finally, we obtain that the estimated number satisfy the following bound
$$
N_h\leqslant n\frac{2n((A-1)+(B-1)t+(C-1)t)}{D}<32n^5h^{3/4}t^{2/3}.\Box
$$

\section{Proof of Theorem \ref{th-energy}}

Let us estimate the number $E_{P}^q(G)$ of solutions $(x_1,y_1,\dots,x_q,y_q)$ of the system
\begin{eqnarray}\label{syst-energy}
P(x_1,y_1)=P(x_2,y_2)=\dots=P(x_q,y_q),\qquad x_i,y_i\in G,\quad i=1,\dots,q.
\end{eqnarray}
The number of solutions of (\ref{syst-energy}) is equal to the following sum
\begin{eqnarray}\label{sum-energy}
E_{P}^q(G)=\sum_{c\in\mathbb{F}_p} (\#\{ (x,y)\mid P(x,y)=c;\,\,\, x,y\in G  \})^q
\end{eqnarray}
The number of solutions $(x,y)$ of equation 
\begin{eqnarray}\label{poly-equ-sred}
P(x,y)=C,
\end{eqnarray}
such that $x,y\in G$ does not exceed $16n^3|G|^{2/3}$.

Consider cosets $r_1G,\dots,r_sG$, $s=\frac{p-1}{t}$.
\begin{eqnarray}\label{poly-syst-est}
P(x,y)=C,\qquad x,y\in G.
\end{eqnarray}
Let us separate equations (\ref{poly-syst-est}) by $|G|$ groups
$$
P(x,y)=gr_i,\qquad x,y\in G,\quad i=1,\ldots,s,
$$
for each $g\in G$.

On the one hand we estimate the sum $\sum_{j=1}^{h} \#\{ (x,y)\mid P(x,y)=hg_{i_j};\,\,\, x,y\in G  \}\leqslant 32h^{3/4}n^5|G|^{2/3}$ by Theorem \ref{th-sr}, and on the other hand the total number of solutions of all equations (\ref{poly-syst-est}) (for all $C\in\mathbb{F}_p$) does not exceed $|G|^2$. It means that the following holds
$$
\sum_{j=1}^{h} (\#\{ (x,y)\mid P(x,y)=gr_{i_j};\,\,\, x,y\in G  \})^q\leqslant \frac{3^q\cdot 2^{3q+2}n^{5q}}{4-q}h^{1-q/4}|G|^{2q/3},\qquad q\leqslant 3.
$$
is the Let us consider the case when this sum the largest. It is easy to see that this case is reached if $h=\frac{1}{32n^4}t^{1/3}$ for each of $t$ cosets.
Now we obtain the following upper bounds for $E_P^q(G)$:
$$
E_P^q(G)\leqslant |G|\frac{3^q\cdot 2^{3q+2}n^{5q}}{4-q}\left(\frac{|G|}{32n^4}\right)^{1-q/4}|G|^{2q/3}<C_1(n,q)|G|^{\frac{7q+16}{12}},\quad q\leqslant 3,
$$
and $C_1(n,q)=\frac{3^q 2^{\frac{17q}{4}-3}n^{6q-4}}{4-q}$. If $q\geqslant 5$, then it is easy to see that $E_P^q(G)<C_2(n,q)|G|^{1+\frac{2q}{3}}$, where $C_2(n,q)=\frac{3^q\cdot 2^{3q+2}n^{5q}}{q-4}$, and  $E_P^4\leqslant C_3(n,q)|G|^{1+\frac{2q}{3}}\ln |G|$, where $C_3(n,q)=3^{q-1}\cdot 2^{3q}n^{5q}$. 
$\Box$

\section{Acknowledgments}

The authors are grateful to Sergey Konyagin, Ilya Shkredov and Ivan Yakovlev for their attention and useful comments. The authors are particularly grateful to Igor Shparlinski and Umberto Zannier for their contribution to the formulation of the problem, which is considered in the paper.

Vyugin I.V.\\
Insitute for Information Transmission Problems RAS,  \\
and\\
National Research University Higher School of Economics,\\
{\it vyugin@gmail.com}.\\
\\
Makarychev S.V.\\
National Research University Higher School of Economics,\\
{\it sergei-lenin2008@yandex.ru}.

\end{document}